\documentclass[twoside,11pt]{article}

\usepackage{blindtext}
\usepackage{subcaption}

\usepackage[utf8]{inputenc} 
\usepackage[T1]{fontenc}    
\usepackage{booktabs}       
\usepackage{amsmath}
\usepackage{amsfonts}       
\usepackage{algpseudocode}
\usepackage{algorithm}
\usepackage{nicefrac}       
\usepackage[abbrvbib,preprint]{jmlr2e}

\usepackage{cleveref}

\usepackage{localmacros}

\usepackage{lastpage}
\jmlrheading{XX}{2025}{1-\pageref{LastPage}}{2/25; Revised X/XX}{X/XX}{21-0000}{Joyce Chew, Willem Diepeveen, and Deanna Needell}

\ShortHeadings{Curvature Corrected Nonnegative Manifold Data Factorization}{Chew, Diepeveen, and Needell}
\firstpageno{1}

\begin{document}

\title{Curvature Corrected Nonnegative Manifold Data Factorization}

\author{\name Joyce Chew \email joycechew@math.ucla.edu \\
       \addr Department of Mathematics\\
       University of California, Los Angeles
       \AND
       \name Willem Diepeveen\email wdiepeveen@math.ucla.edu \\
       \addr Department of Mathematics\\
       University of California, Los Angeles
       \AND
       \name Deanna Needell\email deanna@math.ucla.edu \\
       \addr Department of Mathematics\\
       University of California, Los Angeles}

\editor{}

\maketitle

\begin{abstract}
Data with underlying nonlinear structure are collected across numerous application domains, necessitating new data processing and analysis methods adapted to nonlinear domain structure. Riemannanian manifolds present a rich environment in which to develop such tools, as manifold-valued data arise in a variety of scientific settings, and Riemannian geometry provides a solid theoretical grounding for geometric data analysis. Low-rank approximations, such as nonnegative matrix factorization (NMF), are the foundation of many Euclidean data analysis methods, so adaptations of these factorizations for manifold-valued data are important building blocks for further development of manifold data analysis. In this work, we propose curvature corrected nonnegative manifold data factorization (CC-NMDF) as a geometry-aware method for extracting interpretable factors from manifold-valued data, analogous to nonnegative matrix factorization. We develop an efficient iterative algorithm for computing CC-NMDF and demonstrate our method on real-world diffusion tensor magnetic resonance imaging data.
\end{abstract}


\section{Introduction.}
\label{sec:introduction}
Many scientific techniques collect data with natural representations in non-Euclidean domains. For example, scientific imaging methods such as electron backscatter diffraction imaging \citep{adams1993orientation} and diffusion tensor magnetic resonance imaging \citep{basser1994mr} generate manifold-valued data. For the former, the data live in $SO(3)$, the 3D rotation group, and for the latter, the data live in $\mathcal{P}(3)$, the space of $3 \times 3$ symmetric positive definite matrices. Molecules \citep{kearnes2016molecular}, metabolic pathways \citep{jeong2000large}, and social relationships \citep{newman2002random} are just a few examples of phenomena that can be encoded as graphs. Data collected from a variety of domains can be embedded into Riemannian manifolds to uncover latent non-Euclidean structure \citep{wilson2014spherical}. Additionally, high-dimensional Euclidean data can be modeled as instead being drawn from an unknown, lower-dimensional non-Euclidean domain, enabling new methods of visualization, compression, and further analysis \citep{bronstein2017geometric,puchkin2022structure}. 

In particular, the data setting of \emph{symmetric Riemannian manifolds} is an exciting realm in which to develop data analysis methods for two main reasons. First, the class of these symmetric spaces includes manifolds that appear in many existing scientific techniques, such as those mentioned above,  and manifolds that are popular targets for embeddings of high-dimensional Euclidean data, such as spherical and hyperbolic spaces \citep{lopez2021symmetric,keller2020hydra,wilson2014spherical}. Hence, methods developed for such manifold-valued data and embeddings can be used in a variety of application domains. Second, symmetric spaces have mathematical properties that enable the design of tractable algorithms that respect the geometry of the spaces. Therefore, practical and effective data analysis methods can be developed for these settings. 

For data that have been collected from or embedded in Riemannian manifolds, data analysis methods should account for and leverage the underlying geometry of these manifolds. At the same time, evaluating manifold-specific mappings can be computationally expensive. One approach is to exploit nice geometric properties of a specific manifold domain, such as the sphere \citep{li1999multiscale, narcowich1996nonstationary,tabaghi2023principal}, to enable computationally feasible yet geometrically faithful adaptations of Euclidean data processing techniques. Other methods utilize a \emph{tangent space} of the manifold to approximate the underlying geometry while enjoying the efficiency enabled by working on linear spaces \citep{fletcher2004principal,fletcher2004principal2,ho2013nonlinear,yin2015nonlinear}. 

A foundational family of  techniques for processing Euclidean data is that of \emph{low rank approximation}, motivated by applications such as visualization, compression, feature extraction, and improving the performance of subsequent machine learning methods. Consequently, there is growing interest in developing methods for manifold-valued data that emulate the behavior of low rank approximations and are adapted for these domains. For example, several ideas such as Principal Geodesic Analysis (PGA) \citep{fletcher2004principal,fletcher2004principal2} and the nested subspaces framework \citep{jung2012analysis,yang2021nested} are motivated by the utility of Principal Component Analysis (PCA) for Euclidean data. 

A popular low-rank approximation technique for interpretable feature extraction from nonnegative Euclidean data is \emph{nonnegative matrix factorization} (NMF). For a nonnegative data matrix $\Matrix \in \Real_{\geq 0}^{\numData \times \dimInd}$ and given rank $\numFactors$, NMF is an approximation
\begin{equation}
    \Matrix \approx \mathbf{H} \mathbf{W} 
\end{equation}
where $\mathbf{H} \in \Real_{\geq 0}^{\numData \times \numFactors}, \mathbf{W} \in \Real_{\geq 0}^{\numFactors \times \dimInd}$. By enforcing a nonnegativity constraint on $\mathbf{H}$ and $\mathbf{W}$, NMF extracts additive factors (the rows of $\mathbf{W}$) from high-dimensional Euclidean data that tend to be interpretable by human experts across a variety of data domains, which makes NMF a simple yet powerful tool for tasks such as topic modeling and image segmentation \citep{Lee1999}. One reason for developing methods adapted for manifold-valued data is to retain the structure of the original data domain, so readily interpretable low-rank approximations on manifolds, such as adaptations of NMF, are an exciting avenue for practical applications.

As previously outlined, there are two main types of methods for manifold-valued data. Typically, methods that account for global geometry are specific to particular manifolds in order to leverage known geometric information. While these methods remain faithful to the nonlinear structure of their domains, they are difficult to practically generalize to broader families of manifolds. On the other hand, methods that work on tangent spaces can more easily adapt frameworks originally defined for Euclidean data for many manifold domains. However, tangent space-based methods generally introduce large errors, since curvature information is lost when the data are linearized. These drawbacks present considerable challenges to developing a general NMF-type factorization for the manifold setting. 

In this work, we propose \emph{curvature corrected nonnegative manifold data factorization} (CC-NMDF) as a geometry-aware analogue of NMF for manifold-valued data. We formulate an interpretable factorization for data drawn from any symmetric Riemannian manifold, apply a curvature correction to a fast tangent space-based algorithm, and propose an iterative algorithm for computing CC-NMDF. We demonstrate our method on a real-world data set collected via diffusion tensor magnetic resonance imaging.

\subsection{Related Work.}
In this section, we survey relevant connections between low-rank approximation methods for manifold-valued data and NMF. As previously mentioned, PGA \citep{fletcher2004principal, fletcher2004principal2} is a general method that adapts PCA to the manifold setting, and practical algorithms for implementing PGA reduce to performing PCA on a particular tangent space of the underlying manifold. The learned geodesic submanifolds are intended to maximize reconstruction fidelity, analogous to the role of the linear subspaces learned by PCA. Other approaches construct nested submanifolds using manifold-specific properties of settings such as spherical and hyperbolic surfaces \citep{jung2012analysis,tabaghi2023principal} and Grassmann manifolds \citep{curry2019principal,yang2021nested} to achieve computationally feasible algorithms. However, these methods do not yield submanifolds that are easily interpreted in the original data domain. 

In the NMF literature, some methods have been proposed to account for possible nonlinear structure when computing otherwise standard Euclidean NMF. For instance, \cite{cai2008non}, \cite{huang2014robust}, and \cite{lu2012manifold} construct a nearest-neighbors graph using affinities between input data points (which are high-dimensional vectors in $\Real^D$) and use the corresponding graph Laplacian to regularize NMF in hopes of preserving unknown manifold structure in the data. \cite{he2020low} propose a similarly regularized NMF whose basis matrix is additionally assumed to lie on a Stiefel manifold in order to reduce redundancy in the basis factors. These methods do not have prior knowledge of the underlying manifold the data may be drawn from; instead, they seek to approximate the unknown  geometry through data-driven graphs. Hence, the resulting factors cannot be guaranteed to lie on the underlying manifold, and these methods do not incorporate known geometry of application domains. 

There has been some work on NMF-type methods for data drawn from known manifolds. For example, \cite{ho2013nonlinear} propose a version of sparse coding and dictionary learning for manifold-valued data. However, their methods enforce sparsity rather than nonnegativity for the learned coefficients, so while the dictionary atoms do lie on the manifold the data are drawn from, they cannot be easily interpreted themselves. Also, the algorithms provided require computation of geodesics, so feasible algorithms are provided for only two manifold settings. 

\subsection{Contributions.}
In this work, we adapt a framework for curvature corrected low-rank approximation of manifold-valued data \citep{diepeveen2023curvature} to develop a geometry-aware yet tangent space-based factorization for data drawn from symmetric Riemannian manifolds that is analogous to NMF. Then, we propose an interpretation for the resulting factors in the original manifold domain. We develop an iterative algorithm to compute our CC-NMDF that only requires a multiplicative update and solving a linear system. We also discuss adjustments to the algorithm to improve the interpretability of the computed factors. We demonstrate our method on a real-world dataset and compare it with other low-rank approximation methods to investigate curvature effects and the quality of the resulting factors. 

\subsection{Outline.}
In \Cref{sec:notation} we summarize notation used throughout this paper and necessary concepts from Riemannian geometry. Then, in \Cref{sec:factorization} we introduce a formulation for nonnegative factorization of manifold-valued data, and we apply a curvature correction to account for the underlying geometry of the domain, resulting in our CC-NMDF. Next, in \Cref{sec:algorithms}, we provide tractable algorithms for computing CC-NMDF for data drawn from any symmetric Riemannian manifold, and we discuss algorithm considerations such as initialization. In \Cref{sec:numerics}, we apply our method to a diffusion tensor magnetic resonance imaging dataset, evaluate the effects of algorithm parameter choices, and compare CC-NMDF to other low-rank approximation methods for manifold-valued data on the basis of reconstruction and interpretability of the resulting factors. Finally, in \Cref{sec:conclusions}, we summarize our findings.

\section{Notation and Preliminaries.}
\label{sec:notation}
We first introduce necessary concepts from Riemannian geometry; for more detailed treatment and further reading, see \cite{boothby2003introduction,carmo1992riemannian,lee2013smooth,sakai1996riemannian}. 

Throughout this work, we let $\manifold$ denote a $\dimInd$-dimensional symmetric \emph{Riemannian manifold}. The \emph{tangent space} at a point $\mPoint \in \manifold$ is denoted by $\tangent_\mPoint \manifold$, and it is the space of all derivations at $\mPoint$. The \emph{inner product} or \emph{metric tensor} at a point $\mPoint$ is denoted by $(\cdot, \cdot)_\mPoint: \tangent_\mPoint \manifold \times \tangent_\mPoint \manifold \mapsto \Real$ and defines a Riemannian manifold $(\manifold, (\cdot, \cdot))$, and we denote the associated \emph{Riemannian metric} by $\distance_\manifold: \manifold \times \manifold \to \Real_{\geq 0}$. We denote the \emph{geodesic}, or curve of minimal length, between $\mPoint, \mPointB \in \manifold$ by $\geodesic_{\mPoint, \mPointB}:[0,1] \to \manifold$. For a given point $\mPoint \in \manifold$ and tangent vector $\Xi_\mPoint \in \tangent_\mPoint \manifold$, we define the curve $t \to \geodesic_{\mPoint, \Xi_\mPoint} (t)$ to be the geodesic starting at $\mPoint$ with velocity $\dot{\geodesic}_{\mPoint, \Xi_\mPoint} = \Xi_\mPoint$. These geodesics are used to define the \emph{exponential map} $\exp_\mPoint : \mathcal{G}_\mPoint \manifold \to \manifold$ where
\begin{equation}
    \exp_\mPoint (\Xi_\mPoint) \coloneqq \geodesic_{\mPoint, \Xi_\mPoint}(1)
\end{equation}
and $\mathcal{G}_\mPoint \subset \tangent_\mPoint \manifold$ is
\begin{equation}
    \mathcal{G}_\mPoint \coloneqq \{\Xi_\mPoint: \geodesic_{\mPoint, \Xi_\mPoint}(1) \text{ is defined}\}.
\end{equation}
We will consider \emph{geodesically complete} manifolds throughout, which satisfy $\mathcal{G}_\mPoint = \tangent_\mPoint \manifold$. The inverse of $\exp_\mPoint$ is the \emph{logarithmic map} $\log_\mPoint : \exp_\mPoint (\mathcal{G}_\mPoint') \to \mathcal{G}_\mPoint$, where $\mathcal{G}_\mPoint'\subset \mathcal{G}_\mPoint$ is the set on which $\exp_\mPoint$ is a diffeomorphism. Since $\log_\mPoint$ maps from the manifold $\manifold$, which has nonlinear structure in general, to the tangent space $\tangent_\mPoint \manifold$, which is a vector space, application of the logarithmic map can be considered a linearization. Finally, throughout this work we consider a dataset $\Tensor=\{\Tensor^\sumIndA\}_{\sumIndA=1}^\numData$ where each $\Tensor^\sumIndA \in \manifold$.

\section{Interpretable Factorization of Manifold-Valued Data.}
\label{sec:factorization}
We first seek to adapt the Euclidean NMF to the manifold setting in hopes of obtaining similarly interpretable factors. However, there are two main obstacles to this goal. First, the original Euclidean NMF is only defined for nonnegative input data, and there is no characteristic clearly equivalent to nonnegativity for general manifold-valued data. Second, the interpretability of NMF is derived from its additive nature, which depends on the linearity of the underlying domain, and the non-Euclidean settings that require the use of tools from Riemannian geometry are clearly nonlinear. We address the first obstacle by turning to an NMF variant known as semi-nonnegative matrix factorization \citep{ding2008convex}, which we refer to as semi-NMF, and is defined for (Euclidean) data of mixed sign. We address the second obstacle by proposing interpretation of resulting manifold-valued factors that does not rely on linear structure in the data domain. 
\subsection{Euclidean Semi-NMF and Exact Formulation on Manifolds.}
For Euclidean data of mixed sign, \cite{ding2008convex} proposed a semi-nonnegative matrix factorization that only constrains the coefficients to be nonnegative and allows the learned basis factors to be of mixed sign. Due to the nonnegativity of the coefficients, the resulting factorizations tend to retain the interpretability of standard NMF. Given a data matrix $\Matrix \in \Real^{\numData \times \dimInd}$ and a desired rank $\numFactors$, semi-NMF computes the low-rank approximation
\begin{equation}
\label{eq:eucl SNMF}
    \Matrix \approx \nonnegFactor \factorCoords
\end{equation}
where $\nonnegFactor \in \Real_{\geq 0}^{\numData \times \numFactors}$ and $\factorCoords\in\Real^{\numFactors \times \dimInd}$ and $\nonnegFactor$ and $\factorCoords$ are given by

\begin{equation}
\label{eq:eucl SNMF opt prob}
\underset{\substack{\nonnegFactor \in \Real_{\geq 0}^{\numData \times \numFactors}\\ \factorCoords\in\Real^{\numFactors \times \dimInd}}}{\argmin} \left\|\Matrix - \nonnegFactor \factorCoords\right\|_F^2.
\end{equation}

If we consider $\Matrix$ as a collection of $\numData$ points in $\Real^\dimInd$, we can use the idea of only constraining the learned coefficients for each point to be nonnegative to formulate a similar factorization of a manifold-valued dataset $\{\Tensor^\sumIndA\}_{\sumIndA = 1}^\numData \subset \manifold$ as follows. Given a desired rank $\numFactors$ and a point of linearization $\mPoint \in \manifold$, we seek to find

\begin{equation}
    \label{eq:manifold SNMF exact}
    \underset{\substack{\nonnegFactor \in \Real_{\geq 0}^{\numData \times \numFactors} \\
    \mTVectorFactor_\mPoint^1, \dots, \mTVectorFactor_\mPoint^\numFactors \in \tangent_\mPoint\manifold}}{\argmin} \sum_{i=1}^\numData \distance_\manifold \left(\Tensor^\sumIndA, \exp_\mPoint \left(\sum_{\sumIndC=1}^\numFactors \nonnegFactor_{\sumIndA,\sumIndC} \mTVectorFactor_\mPoint^\sumIndC\right)\right)^2.
\end{equation}

Defining such a factorization in this way requires that the learned factors actually lie on a chosen tangent space $\tangent_\mPoint \manifold$. Note that this formulation reduces to the Euclidean SNMF when $\manifold$ is $\Real^\dimInd$ and $\mPoint = 0$. In this case, the learned tangent vector factors $\mTVectorFactor_\mPoint^\sumIndC$ are the rows of $\factorCoords$. 
In the general case, we define the corresponding manifold-valued factors $\{\TensorY^\sumIndC\}_{\sumIndC = 1}^\numFactors$ by
\begin{equation}
\label{eq:manifold valued factors}
    \TensorY^\sumIndC \coloneqq \exp_\mPoint \left ( \left(\max_{\sumIndA = 1, \dots, \numData} \nonnegFactor_{\sumIndA, \sumIndC} \right) \mTVectorFactor_\mPoint^\sumIndC \right )
\end{equation}
and we interpret these manifold-valued factors as the vertices of a nonlinear polyhedron in $\manifold$ that approximately encapsulates the original data. Then, the learned coefficients for each of the original data points, i.e. the rows of $\nonnegFactor$, locate each data point with respect to each vertex and the base point of linearization and hence serve as geometric descriptors. 

While the formulation given in \eqref{eq:manifold SNMF exact} is a natural extension of the Euclidean semi-NMF to the setting of manifold-valued data, in practice it is computationally expensive to solve via direct optimization on $\manifold$. In the following section, we develop feasible approximations to the exact problem.

\subsection{Tangent Space and Curvature Corrected Nonnegative Manifold Data Factorizations.}
\label{sec:tsnmdf}
Because of the computational expense of direct optimization on $\manifold$, we first propose to solve an approximation to the exact problem \eqref{eq:manifold SNMF exact} in $\tangent_\mPoint \manifold$. We rewrite the exact problem as
\begin{equation}
     \label{eq:manifold snmf rewritten}   \argmin_{\substack{\nonnegFactor \in \Real_{\geq 0}^{\numData \times \numFactors}\\
    \mTVectorFactor_\mPoint^1, \dots, \mTVectorFactor_\mPoint^\numFactors \in \tangent_\mPoint \manifold} } \sum_{\sumIndA=1}^\numData  \distance_\manifold\left( \exp_\mPoint\left( \log_\mPoint \Tensor^\sumIndA\right), \exp_\mPoint \left(\sum_{\sumIndC=1}^\numFactors \nonnegFactor_{\sumIndA,\sumIndC} \mTVectorFactor_\mPoint^\sumIndC\right)\right)^2.
\end{equation}
A natural approach to approximating this formulation is to linearize the problem and solve
\begin{equation}
\label{eq:tnmdf problem}
    \argmin_{\substack{\nonnegFactor \in \Real_{\geq 0}^{\numData \times \numFactors}\\
    \mTVectorFactor_\mPoint^1, \dots, \mTVectorFactor_\mPoint^\numFactors \in \tangent_\mPoint \manifold} } \sum_{\sumIndA=1}^\numData   \left \| \log_\mPoint \Tensor^\sumIndA - \sum_{\sumIndC=1}^\numFactors \nonnegFactor_{\sumIndA, \sumIndC} \mTVectorFactor_\mPoint^\sumIndC\right\|_\mPoint^2.
\end{equation}
We refer to the loss function considered in \eqref{eq:tnmdf problem} as the \emph{tangent space error}, as it is the error of approximation in $\tangent_\mPoint \manifold$. Accordingly, we refer to the corresponding approximation 
\begin{equation}
\label{eq:tnmdf}
    \log_\mPoint \Tensor^\sumIndA \approx \sum_{\sumIndC = 1}^\numFactors \nonnegFactor_{\sumIndA, \sumIndC} \mTVectorFactor_\mPoint ^\sumIndC,
\end{equation}
where $\nonnegFactor$ and $\{\mTVectorFactor_\mPoint^\sumIndC\}_{\sumIndC = 1}^\numFactors \subset \tangent_\mPoint \manifold$ are solutions of \eqref{eq:tnmdf problem}, as \emph{tangent space nonnegative manifold data factorization} (T-NMDF), and we define associated manifold-valued factors as in \eqref{eq:manifold valued factors}. Crucially, T-NMDF can be computed solely on $\tangent_\mPoint \manifold$, which greatly speeds up computation due to linearity of $\tangent_\mPoint \manifold$. However, the linearization of the data required to work solely on $\tangent_\mPoint \manifold$ may incur additional approximation error due to the global geometry of $\manifold$ \citep{diepeveen2023curvature}. In particular, the solutions of \eqref{eq:tnmdf problem} are not guaranteed to coincide with the solutions of the exact problem \eqref{eq:manifold SNMF exact}. Hence, T-NMDF may produce reconstructions that are not as faithful to the original data as solutions of \eqref{eq:manifold SNMF exact}, and the manifold-valued factors learned by T-NMDF may suffer distortions that impede their interpretability in the original data domain. To correct these issues while retaining the computational feasibility of T-NMDF, we propose a \emph{curvature corrected nonnegative manifold data factorization} (CC-NMDF) as follows.

By theory developed by \citealp{diepeveen2023curvature}, we know that the exact reconstruction error on $\manifold$ minimized in \eqref{eq:manifold SNMF exact} is dominated by \emph{curvature corrected} tangent space error, rather than the linearized error considered in \eqref{eq:tnmdf problem}. In particular, consider a symmetric Riemannian manifold $(\manifold, (\cdot, \cdot) )$ with \emph{curvature tensor} at a point $\mPoint$ denoted by $\curvature_\mPoint(\cdot, \cdot)(\cdot):\tangent_\mPoint \manifold \times \tangent_\mPoint \manifold \times \tangent_\mPoint \manifold \to \tangent_\mPoint \manifold$. For each $\sumIndA \in [\numData]$ let $\{ \Theta_{\mPoint}^{(\sumIndA),\sumIndB} \}_{\sumIndB=1}^\dimInd \subset \tangent_\mPoint \manifold$ be an orthonormal frame that diagonalizes the operator
\begin{equation}
    \Theta_\mPoint \mapsto \curvature_\mPoint \left(\Theta_\mPoint, \log_{\mPoint} \Tensor^\sumIndA\right) \log_{\mPoint} \Tensor^\sumIndA
\end{equation}
with respective eigenvalues $\kappa_{\sumIndA,\sumIndB}$ and define $\beta: \Real \to \Real$ as
\begin{equation}
\label{eq:beta def}
    \beta(\kappa):=\left\{\begin{array}{cl}
\frac{\sinh (\sqrt{-\kappa})}{\sqrt{-\kappa}}, & \kappa<0, \\
1, & \kappa=0, \\
\frac{\sin (\sqrt{\kappa})}{\sqrt{\kappa}}, & \kappa>0 .
\end{array}\right.
\end{equation}
Then, we have by \cite{diepeveen2023curvature}, Thm.~3.4 that
\begin{equation}
    \sum_{\sumIndA=1}^\numData  \distance_{\manifold} \left(\Tensor^\sumIndA, \exp_{\mPoint}\left(\sum_{\sumIndC=1}^\numFactors \nonnegFactor_{\sumIndA,\sumIndC} \mTVectorFactor_\mPoint^\sumIndC\right)\right)^2 = \sum_{\sumIndA=1}^\numData  \sum_{\sumIndB=1}^\dimInd \beta(\kappa_{\sumIndA,\sumIndB})^2 \left(\sum_{\sumIndC=1}^\numFactors \nonnegFactor_{\sumIndA,\sumIndC} \mTVectorFactor_\mPoint^\sumIndC - \log_{\mPoint} \Tensor^\sumIndA, \Theta_{\mPoint}^{(\sumIndA),\sumIndB}\right)_\mPoint^2 + \mathcal{O}(\epsilon^3)
    \label{eq:relaxed-semi-nmf}
\end{equation}
where 
\begin{equation}
    \epsilon \coloneqq \sqrt{\sum_{\sumIndA=1}^\numData \left\|\sum_{\sumIndC=1}^\numFactors \nonnegFactor_{\sumIndA,\sumIndC} \mTVectorFactor_\mPoint^\sumIndC - \log_{\mPoint} \Tensor^\sumIndA\right\|_\mPoint^2}.
\end{equation}
Following \cite{diepeveen2023curvature}, we refer to the leading term of the right hand side of \eqref{eq:relaxed-semi-nmf} as the \emph{curvature corrected approximation error}. 
This theory implies that when $\manifold$ has nonzero curvature, the curvature corrected approximation error is a better approximation of the exact reconstruction error on $\manifold$ than the tangent space error. In fact, \cite{diepeveen2023curvature} found that for manifold-valued singular value decompositions, minimizing the curvature corrected approximation error resulted in better performance (on the basis of exact reconstruction error) than minimizing the tangent space error.  
This suggests that solving the curvature corrected problem
\begin{equation}
    \label{eq:cc snmdf prob}
    \argmin_{\substack{\nonnegFactor \in \Real_{\geq 0}^{\numData \times \numFactors} \\
        \mTVectorFactor_\mPoint^1, \dots, \mTVectorFactor_\mPoint^\numFactors \in \tangent_\mPoint \manifold}} \sum_{\sumIndA=1}^\numData  \sum_{\sumIndB=1}^\dimInd \beta(\kappa_{\sumIndA,\sumIndB})^2 \left(\sum_{\sumIndC=1}^\numFactors \nonnegFactor_{\sumIndA,\sumIndC} \mTVectorFactor_\mPoint^\sumIndC - \log_{\mPoint} \Tensor^\sumIndA, \Theta_{\mPoint}^{(\sumIndA),\sumIndB}\right)_\mPoint^2
\end{equation}
will give a geometry-aware approximate minimizer of \eqref{eq:manifold SNMF exact} that only requires working on $\tangent_\mPoint \manifold$ rather than $\manifold$ itself. We refer to the corresponding approximation as \emph{curvature corrected nonnegative manifold data factorization} (CC-NMDF). 
\section{Algorithms for Computing Nonnegative Manifold Data Factorizations.}
\label{sec:algorithms}
In this section, we present iterative algorithms for computing T-NMDF and CC-NMDF, and we propose a data-driven initial guess for CC-NMDF. We also discuss postprocessing and parameter choices to improve the interpretability of the learned factorizations.
\subsection{Tangent Space Nonnegative Manifold Data Factorization.}
To find an approximate solution of \cref{eq:tnmdf problem}, we observe that we can apply algorithms for computing Euclidean semi-NMF. Explicitly, we fix an orthonormal basis $\{\phi_\mPoint^\sumIndB\}_{\sumIndB=1}^\dimInd \subset \tangent_\mPoint \manifold$ and consider the coordinate matrix $\Matrix \in \Real^{\numData \times \dimInd}$ defined entrywise by
\begin{equation}
    \Matrix_{\sumIndA, \sumIndB} = \left ( \log_\mPoint \Tensor^\sumIndA, \phi_\mPoint^\sumIndB \right )_\mPoint
\end{equation}
so that
\begin{equation}
\label{eq:coord matrix}
    \log_{\mPoint} \Tensor^\sumIndA = \sum_{\sumIndB=1}^\dimInd \Matrix_{\sumIndA,\sumIndB} \phi_\mPoint^\sumIndB.
\end{equation}
Then, for a desired rank $\numFactors$, we compute the semi-NMF of $\Matrix$ given by
\begin{equation}
    \Matrix \approx \nonnegFactor \factorCoords
\end{equation}
where $\nonnegFactor \in \Real_{\geq 0}^{\numData \times \numFactors}$ and $\factorCoords\in\Real^{\numFactors \times \dimInd}$, using an iterative algorithm such as the one described by \cite{ding2008convex}. Then, the corresponding tangent vector factors are 
$\{\mTVectorFactor_\mPoint^\sumIndC\}_{\sumIndC=1}^\numFactors \subset \tangent_\mPoint\manifold$, where $\mTVectorFactor_\mPoint^\sumIndC = \sum_{\sumIndB=1}^\dimInd \factorCoords_{\sumIndC, \sumIndB}\phi_{\mPoint}^\sumIndB$, and the manifold-valued factors are $\{\TensorY^\sumIndC\}_{\sumIndC=1}^\numFactors \subset \manifold$ where
\begin{equation}
    \TensorY^\sumIndC = \exp_\mPoint \left ( \left(\max_{\sumIndA = 1, \dots, \numData} \nonnegFactor_{\sumIndA, \sumIndC}\right) \mTVectorFactor_\mPoint^\sumIndC \right ).
\end{equation}
The T-NMDF algorithm (using the Euclidean semi-NMF algorithm from \citealp{ding2008convex}) is summarized in \Cref{alg:t-SNMDF}. This approach is a computationally cheap way to find approximate solutions of \eqref{eq:tnmdf problem}, since it requires only standard matrix operations after the coordinate matrix $\Matrix$ is obtained. 

\begin{algorithm}[t!]
\caption{Tangent space nonnegative manifold data factorization (T-NMDF)}
\label{alg:t-SNMDF}
\begin{algorithmic}[1]
\Require{$\{\Tensor^\sumIndA\}_{\sumIndA=1}^\numData \subset \manifold$, $\mPoint\in\manifold$, $K \in \Natural$}
\For{$\sumIndA = 1, \ldots, \numData$}
\For{$\sumIndB = 1, \ldots, \dimInd$}
\State $\Matrix_{\sumIndA, \sumIndB} \gets \left ( \log_\mPoint \Tensor^\sumIndA, \phi_\mPoint^\sumIndB \right )_\mPoint$
\EndFor
\EndFor
\State Compute the Euclidean semi-NMF $\Matrix \approx \nonnegFactor \factorCoords$
\For{$\sumIndC = 1, \ldots, \numFactors$}
\State $\mTVectorFactor_\mPoint^\sumIndC \coloneqq \sum_{\sumIndB=1}^\dimInd \factorCoords_{\sumIndC, \sumIndB}\phi_{\mPoint}^\sumIndB$
\State $\TensorY^\sumIndC \coloneqq \exp_\mPoint \left ( \left(\max_{\sumIndA = 1, \dots, \numData} \nonnegFactor_{\sumIndA, \sumIndC}\right) \mTVectorFactor_\mPoint^\sumIndC \right )$
\EndFor
\State \Return $\nonnegFactor, \{\mTVectorFactor_\mPoint^\sumIndC\}_{\sumIndC=1}^\numFactors$, and $\{\TensorY^\sumIndC\}_{\sumIndC=1}^\numFactors$
\end{algorithmic}
\end{algorithm}

\subsection{Curvature Corrected Nonnegative Manifold Data Factorization.}
\label{sec:ccnmdf}
We propose an alternating iterative algorithm for solving the curvature-aware factorization problem \eqref{eq:cc snmdf prob}. First, we argue that the solution of tangent space $\numFactors$-means will give a good initialization. Then, we derive a multiplicative update for the curvature corrected algorithm.

\subsubsection{Initialization of the algorithm.}
\label{sec:initialization}
The framework developed by \cite{diepeveen2023curvature} implies that good solutions to low-rank decomposition problems such as \eqref{eq:cc snmdf prob} can be obtained by first solving a similar problem in $\tangent_\mPoint \manifold$ and subsequently applying curvature correction. In particular, we claim that CC-NMDF can be obtained as a relaxation of the manifold K-means problem. Consider the $\numFactors$-means clustering problem for data $\{ \Tensor^\sumIndA\}_{\sumIndA = 1}^\numData \subset \manifold$
\begin{equation}
\label{eq:exact kmeans}
    \argmin_{\substack{ \mPointB^1, \dots, \mPointB^\numFactors \in \manifold \\
    \nonnegFactor \in \{0,1\}^{\numData \times \numFactors} \\
    \text{s.t. } \nonnegFactor \mathbf{1}_\numFactors = 1}} \sum_{\sumIndA = 1}^{\numData} \sum_{\sumIndC = 1}^{\numFactors}  \nonnegFactor_{\sumIndA, \sumIndC}\distance_\manifold \left(\Tensor^\sumIndA, \mPointB^\sumIndC\right)^2.
\end{equation}
which assigns each data point $\Tensor^\sumIndA$ to one of $\numFactors$ clusters, whose centroids are given by the minimizers $\{\mPointB^\sumIndC\}_{\sumIndC=1}^\numFactors \subset \manifold$ of the above. Then, if $\nonnegFactor$ and $\{\mPointB^\sumIndC\}_{\sumIndC=1}^\numFactors $ are solutions of the above, we see that 

\begin{multline}
    \sum_{\sumIndA=1}^\numData  \sum_{\sumIndC=1}^\numFactors \nonnegFactor_{\sumIndA, \sumIndC} \distance_{\manifold} \left(\Tensor^\sumIndA, \mPointB^\sumIndC\right)^2 \overset{\mTVectorFactor_\mPoint^\sumIndC \coloneqq \log_{\mPoint}\mPointB^{\sumIndC}}{=} \sum_{\sumIndA=1}^\numData  \sum_{\sumIndC=1}^\numFactors \nonnegFactor_{\sumIndA, \sumIndC} \distance_{\manifold} \left(\Tensor^\sumIndA, \exp_{\mPoint}\left( \mTVectorFactor_\mPoint^\sumIndC\right)\right)^2 \\
    \overset{\text{\cite[Thm~3.4]{diepeveen2023curvature}}}{\approx} \sum_{\sumIndA=1}^\numData \sum_{\sumIndC=1}^\numFactors \nonnegFactor_{\sumIndA, \sumIndC} \sum_{\sumIndB=1}^\dimInd  \beta(\kappa_{\sumIndA,\sumIndB})^2 \left( \mTVectorFactor_\mPoint^\sumIndC - \log_{\mPoint} \Tensor^\sumIndA, \Theta_{\mPoint}^{(\sumIndA),\sumIndB}\right)_\mPoint^2\\
    \overset{\nonnegFactor_{\sumIndA, \sumIndC} \in \{0,1\}}{=} \sum_{\sumIndA=1}^\numData \sum_{\sumIndC=1}^\numFactors \sum_{\sumIndB=1}^\dimInd  \beta(\kappa_{\sumIndA,\sumIndB})^2 \left(\nonnegFactor_{\sumIndA, \sumIndC} \mTVectorFactor_\mPoint^\sumIndC -  \nonnegFactor_{\sumIndA, \sumIndC} \log_{\mPoint} \Tensor^\sumIndA, \Theta_{\mPoint}^{(\sumIndA),\sumIndB}\right)_\mPoint^2\\
    \overset{\sum_{\sumIndC=1}^\numFactors \nonnegFactor_{\sumIndA, \sumIndC}=1}{=} \sum_{\sumIndA=1}^\numData  \sum_{\sumIndB=1}^\dimInd  \beta(\kappa_{\sumIndA,\sumIndB})^2 \left(\sum_{\sumIndC=1}^\numFactors\nonnegFactor_{\sumIndA, \sumIndC} \mTVectorFactor_\mPoint^\sumIndC -  \sum_{\sumIndC=1}^\numFactors\nonnegFactor_{\sumIndA,\sumIndC} \log_{\mPoint} \Tensor^\sumIndA, \Theta_{\mPoint}^{(\sumIndA),\sumIndB}\right)_\mPoint^2 \\
\overset{\sum_{\sumIndC=1}^\numFactors \nonnegFactor_{\sumIndA, \sumIndC}=1}{=} \sum_{\sumIndA=1}^\numData  \sum_{\sumIndB=1}^\dimInd  \beta(\kappa_{\sumIndA,\sumIndB})^2 \left(\sum_{\sumIndC=1}^\numFactors\nonnegFactor_{\sumIndA, \sumIndC} \mTVectorFactor_\mPoint^\sumIndC -  \log_{\mPoint} \Tensor^\sumIndA, \Theta_{\mPoint}^{(\sumIndA),\sumIndB}\right)_\mPoint^2 .
\end{multline}
By relaxing the constraint on $\nonnegFactor$ to allow its entries to take on values in $(0, \infty)$ and removing the row sum constraint, we then obtain the optimization problem for CC-NMDF \eqref{eq:cc snmdf prob}. Hence, CC-NMDF can be viewed as a relaxation of the exact $\numFactors$-means problem on $\manifold$, which implies that a solution of manifold $\numFactors$-means is a good initialization for solving the CC-NMDF problem. To avoid the expense of direct computation on $\manifold$ for initialization, we can instead solve the tangent space $\numFactors$-means problem 
\begin{align}
     \label{eq:kmeans tpm naive}   \argmin_{\substack{\nonnegFactor \in \{0, 1\}^{\numData \times \numFactors}\\
    \mTVectorFactor_\mPoint^1, \dots, \mTVectorFactor_\mPoint^\numFactors \in \tangent_\mPoint \manifold} } & \sum_{\sumIndA=1}^\numData \sum_{\sumIndC=1}^\numFactors \nonnegFactor_{\sumIndA, \sumIndC} \left\| \log_\mPoint \Tensor^\sumIndA - \mTVectorFactor_\mPoint^\sumIndC\right\|_\mPoint^2 \\
    \text{subject to }& \nonnegFactor \mathbf{1}_\numFactors = \mathbf{1}_\numData \nonumber
\end{align}
by applying standard Euclidean $\numFactors$-means algorithms to the coordinate matrix $\Matrix$, computed as in \eqref{eq:coord matrix}. Once we obtain an initial guess for $\nonnegFactor$ by tangent space $\numFactors$-means, we relax $\nonnegFactor$ by replacing each 0 entry with $\delta$ for some $0 < \delta \ll 1$, then normalizing the rows of $\nonnegFactor$ such that $\nonnegFactor \mathbf{1}_\numFactors = \mathbf{1}_\numData$.
\subsubsection{Multiplicative Update for CC-NMDF.}
Starting from initial guesses for $\nonnegFactor \in \Real_{\geq 0}^{\numData \times \numFactors}$ and $\{\mTVectorFactor_\mPoint^\sumIndC\}_{\sumIndC = 1}^\numFactors$, we apply curvature correction to  $\nonnegFactor$ and each $\mTVectorFactor_\mPoint^\sumIndC$ in an alternating fashion. First, we correct the tangent vector factors $\mTVectorFactor_\mPoint^\sumIndC$ by solving

\begin{equation}\label{eq:cc-factor-problem}
    \argmin_{\mTVectorFactor_\mPoint^1, \ldots, \mTVectorFactor_\mPoint^\numFactors \in \tangent_\mPoint \manifold} \sum_{\sumIndA=1}^\numData  \sum_{\sumIndB=1}^\dimInd \beta(\kappa_{\sumIndA,\sumIndB})^2 \left(\sum_{\sumIndC=1}^\numFactors \nonnegFactor_{\sumIndA, \sumIndC} \mTVectorFactor_\mPoint^\sumIndC - \log_{\mPoint} \Tensor^\sumIndA, \Theta_{\mPoint}^{(\sumIndA),\sumIndB}\right)_\mPoint^2.
\end{equation}
To do this, we first expand each $\mTVectorFactor_\mPoint^\sumIndC$ and write
\begin{equation}
    \mTVectorFactor_\mPoint^\sumIndC =  \sum_{\sumIndB=1}^\dimInd \factorCoords_{\sumIndC, \sumIndB} \phi_\mPoint^\sumIndB
\end{equation}
where $\factorCoords \in \Real^{\numFactors \times \dimInd}$ is a real-valued matrix and $\{\phi_\mPoint^\sumIndB\}_{\sumIndB=1}^\dimInd \subset \tangent_\mPoint \manifold$ is an orthonormal basis. Then we can solve \eqref{eq:cc-factor-problem} by adapting the curvature correction step in \cite[Sec.~5.2]{diepeveen2023curvature}. Explicitly, we define the real-valued tensor $\TensorB \in \Real^{\numData \times \dimInd \times \numFactors \times \dimInd}$ entrywise by

\begin{equation}
    \label{eq:b tensor def}
    \TensorB_{\sumIndA_1, \sumIndB, \sumIndC_1, \sumIndA} \coloneqq \nonnegFactor_{\sumIndA_1, \sumIndC_1} \left(\phi_\mPoint ^\dimInd, \Theta_{\mPoint}^{(\sumIndA),\sumIndB}\right)_\mPoint
\end{equation}
and then define the real-valued tensor $\TensorA \in \Real^{\numFactors \times \dimInd \times \numFactors \times \dimInd}$ entrywise by

\begin{equation}
    \label{eq:a tensor def}
    \TensorA_{\sumIndC_1', \sumIndA', \sumIndC_1, \sumIndA} \coloneqq \sum_{\sumIndA_1 = 1}^\numData \sum_{\sumIndB = 1}^\dimInd \beta(\kappa_{\sumIndA_1, \sumIndB})^2 \TensorB_{\sumIndA_1, \sumIndB, \sumIndC_1, \sumIndA} \TensorB_{\sumIndA_1, \sumIndB, \sumIndC_1', \sumIndA'}.
\end{equation}
Then, solutions of the linear system
\begin{equation}
    \label{eq:f linear system}
    \sum_{\sumIndC_1 = 1}^\numFactors \sum_{\sumIndA = 1}^\dimInd \TensorA_{\sumIndC_1', \sumIndA', \sumIndC_1, \sumIndA} \factorCoords_{\sumIndC_1, \sumIndA} = \sum_{\sumIndA_1 = 1}^\numData \sum_{\sumIndB = 1}^\dimInd \beta(\kappa_{\sumIndA_1, \sumIndB})^2 \TensorB_{\sumIndA_1, \sumIndB, \sumIndC_1, \sumIndA} \left( \log_\mPoint \Tensor^\sumIndA, \Theta_{\mPoint}^{(\sumIndA),\sumIndB} \right)_\mPoint
\end{equation}
satisfy the first-order optimality conditions for \eqref{eq:cc-factor-problem}.

Next, to correct the non-negative factor $\nonnegFactor$, we seek to solve 
\begin{equation}\label{eq:cc-nonneg-problem}
    \inf_{\nonnegFactor \in \Real_{\geq 0}^{\numData \times \numFactors}} \sum_{\sumIndA=1}^\numData  \sum_{\sumIndB=1}^\dimInd \beta(\kappa_{\sumIndA,\sumIndB})^2 \left(\sum_{\sumIndC=1}^\numFactors \nonnegFactor_{\sumIndA, \sumIndC} \mTVectorFactor_\mPoint^\sumIndC - \log_{\mPoint} \Tensor^\sumIndA, \Theta_{\mPoint}^{(\sumIndA),\sumIndB}\right)_\mPoint^2.
\end{equation}
We derive a multiplicative update for this constrained optimization problem using the Karush-Kuhn-Tucker conditions. Let $f:\Real^{\numData \times \numFactors} \to \Real$ be given by
\begin{equation}
    f(\nonnegFactor) = \sum_{\sumIndA=1}^\numData  \sum_{\sumIndB=1}^\dimInd \beta(\kappa_{\sumIndA,\sumIndB})^2 \left(\sum_{\sumIndC=1}^\numFactors \nonnegFactor_{\sumIndA, \sumIndC} \mTVectorFactor_\mPoint^\sumIndC - \log_{\mPoint} \Tensor^\sumIndA, \Theta_{\mPoint}^{(\sumIndA),\sumIndB}\right)_\mPoint^2.
\end{equation}
We then form the Lagrangian function
\begin{equation}
    g(\nonnegFactor) = -\trace(\lambda \nonnegFactor^T) + f(\nonnegFactor)
\end{equation}
where $\lambda \in \Real^{\numData \times \numFactors}$, enforcing nonnegativity for each entry of $\nonnegFactor$. The gradient of the Lagrangian is given by
\begin{equation}
    \nabla g(\nonnegFactor)_{\sumIndA_1 \sumIndC_1} = -\lambda_{\sumIndA_1 \sumIndC_1} + 2\sum_{\sumIndA=1}^\numData \sum_{\sumIndB=1}^\dimInd \beta(\kappa_{\sumIndA,\sumIndB})^2\left(\sum_{\sumIndC=1}^\numFactors \nonnegFactor_{\sumIndA, \sumIndC} \mTVectorFactor_\mPoint^\sumIndC - \log_{\mPoint} \Tensor^\sumIndA, \Theta_{\mPoint}^{(\sumIndA),\sumIndB}\right)_\mPoint \left(\mTVectorFactor_\mPoint^{\sumIndC_1}, \Theta_{\mPoint}^{(\sumIndA_1),\sumIndB}\right)_\mPoint.
\end{equation}
By the complementary slackness condition, we have that $\lambda_{\sumIndA_1 \sumIndC_1} \nonnegFactor_{\sumIndA_1 \sumIndC_1}=0$ for $\sumIndA_1 = 1, \dots, \numData$ and $\sumIndC_1 = 1, \dots, \numFactors$. Setting $\nabla g(\nonnegFactor)_{\sumIndA_1 \sumIndC_1} = 0$ and using this complementary slackness condition, we obtain the $\numData\numFactors$ equations
\begin{equation}
    2 \nonnegFactor_{\sumIndA_1 \sumIndC_1} \sum_{\sumIndA=1}^\numData \sum_{\sumIndB=1}^\dimInd \beta(\kappa_{\sumIndA,\sumIndB})^2\left(\sum_{\sumIndC=1}^\numFactors \nonnegFactor_{\sumIndA, \sumIndC} \mTVectorFactor_\mPoint^\sumIndC - \log_{\mPoint} \Tensor^\sumIndA, \Theta_{\mPoint}^{(\sumIndA),\sumIndB}\right)_\mPoint \left(\mTVectorFactor_\mPoint^{\sumIndC_1}, (\Theta_{\mPoint}^{(\sumIndA_1),\sumIndB}\right)_\mPoint = 0
\end{equation}
which we equivalently write as the fixed point equations
\begin{equation}
\nonnegFactor_{\sumIndA_1 \sumIndC_1} = \nonnegFactor_{\sumIndA_1 \sumIndC_1} \sqrt{\frac{\mathbf{A}_{\sumIndA_1\sumIndC_1}^+ + \mathbf{B}_{\sumIndA_1\sumIndC_1}^-}{\mathbf{A}_{\sumIndA_1\sumIndC_1}^- + \mathbf{B}_{\sumIndA_1\sumIndC_1}^+}}
\end{equation}
where the entries of $\mathbf{A}$ and $\mathbf{B}$ are given by
\begin{equation}
\label{eq:mult update A}
    \mathbf{A}_{\sumIndA_1, \sumIndC_1} = \sum_{\sumIndA=1}^\numData \sum_{\sumIndB=1}^\dimInd \beta(\kappa_{\sumIndA,\sumIndB})^2\left(\sum_{\sumIndC=1}^\numFactors \nonnegFactor_{\sumIndA,\sumIndC} \mTVectorFactor_\mPoint^\sumIndC, \Theta_{\mPoint}^{(\sumIndA),\sumIndB}\right)_\mPoint \left(\mTVectorFactor_\mPoint^{\sumIndC_1}, (\Theta_{\mPoint}^{(\sumIndA_1),\sumIndB}\right)_\mPoint,
\end{equation}
\begin{equation}
\label{eq:mult update B}
     \mathbf{B}_{\sumIndA_1, \sumIndC_1} = \sum_{\sumIndA=1}^\numData \sum_{\sumIndB=1}^\dimInd \beta(\kappa_{\sumIndA,\sumIndB})^2\left(\log_\mPoint\Tensor^\sumIndA, \Theta_{\mPoint}^{(\sumIndA),\sumIndB}\right)_\mPoint \left(\mTVectorFactor_\mPoint^{\sumIndC_1}, (\Theta_{\mPoint}^{(\sumIndA_1),\sumIndB}\right)_\mPoint,
\end{equation}
and $\mathbf{A}_{\sumIndA_1, \sumIndC_1}^+ = \max\{0, \mathbf{A}_{\sumIndA_1, \sumIndC_1}\}$ and $\mathbf{A}_{\sumIndA_1, \sumIndC_1}^- = -\min\{0, \mathbf{A}_{\sumIndA_1, \sumIndC_1}\}$. Hence, we obtain the multiplicative update 
\begin{equation}
    \nonnegFactor \leftarrow \nonnegFactor \odot \sqrt{\frac{\mathbf{A}^+ + \mathbf{B}^-}{\mathbf{A}^- + \mathbf{B}^+}}
\end{equation}
where $\odot$ denotes pointwise multiplication, and the square root and division are also pointwise. To improve convergence, we repeat the update for $\nonnegFactor$ before re-updating $\factorCoords$. The full algorithm is summarized in \Cref{alg:CC-NMDF}. We note that $\mathbf{B}$ can be pre-computed and used throughout the updates for $\nonnegFactor$, and the computations of $\mathbf{A}$ and $\mathbf{B}$ are both entirely in $\tangent_\mPoint \manifold$.

\begin{algorithm}[t!]
\caption{Curvature corrected nonnegative manifold data factorization (CC-NMDF)}
\label{alg:CC-NMDF}
\begin{algorithmic}[1]
\Require{$\{\Tensor^\sumIndA\}_{\sumIndA=1}^\numData \subset \manifold$, $\mPoint\in\manifold$, $K \in \Natural, \delta$, \texttt{maxIter}, \texttt{maxSubIter}}
\For{$\sumIndA = 1, \ldots, \numData$}
\For{$\sumIndB = 1, \ldots, \dimInd$}
\State $\Matrix_{\sumIndA, \sumIndB} \gets \left ( \log_\mPoint \Tensor^\sumIndA, \phi_\mPoint^\sumIndB \right )_\mPoint$
\EndFor
\EndFor
\State Compute the Euclidean $\numFactors$-means clustering $\Matrix \approx \nonnegFactor \factorCoords$
\State Replace every 0 entry of $\nonnegFactor$ with $\delta$.
\State Set $\mathbf{D} = \text{diag}(\nonnegFactor \mathbf{1}_K)$,  $\nonnegFactor \gets \mathbf{D}^{-1}\nonnegFactor$
\For{$i=1$, \ldots, \texttt{maxIter}}
\State Construct linear system \eqref{eq:f linear system} and update $\factorCoords$ accordingly.
\For{$j=1$, \ldots, \texttt{maxSubIter}}
\State Compute $\mathbf{A}, \mathbf{B}$ as defined in \cref{eq:mult update A,eq:mult update B}.
\State Set $\nonnegFactor \gets \nonnegFactor \odot \nonnegFactor \sqrt{\frac{\mathbf{A}^+ + \mathbf{B}^-}{\mathbf{A}^- + \mathbf{B}^+}}$
\EndFor
\EndFor
\For{$\sumIndC = 1, \ldots, \numFactors$}
\State $\mTVectorFactor_\mPoint^\sumIndC \coloneqq \sum_{\sumIndB=1}^\dimInd \factorCoords_{\sumIndC, \sumIndB}\phi_{\mPoint}^\sumIndB$
\State $\TensorY^\sumIndC \coloneqq \exp_\mPoint \left ( \left(\max_{\sumIndA = 1, \dots, \numData} \nonnegFactor_{\sumIndA, \sumIndC}\right) \mTVectorFactor_\mPoint^\sumIndC \right )$
\EndFor
\State \Return $\nonnegFactor, \{\mTVectorFactor_\mPoint^\sumIndC\}_{\sumIndC=1}^\numFactors$, and $\{\TensorY^\sumIndC\}_{\sumIndC=1}^\numFactors$
\end{algorithmic}
\end{algorithm}

\subsection{Algorithm Considerations.}
\label{sec:algorithm-considerations}
In this section, we discuss two aspects of CC-NMDF that must be considered in practical applications. 
\subsubsection{Potential Cancellation of Tangent Vector Factors.} Because CC-NMDF does not impose restrictions on the learned tangent vector factors $\{\mTVectorFactor_\mPoint^\sumIndC\}_{\sumIndC=1}^\numFactors$, it is possible that cancellation can occur in the learned factors. That is, we may have $(\mTVectorFactor_\mPoint^\sumIndB, \mTVectorFactor_\mPoint^\sumIndC)_\mPoint < 0$ for some $\sumIndB \neq \sumIndC$. This may impact the interpretability of the factors. 

To build intuition about this effect, consider the Euclidean semi-NMF of real-valued data $\Matrix \in \Real^{\numData \times \dimInd}$ with mixed signs 
\begin{equation}
    \label{eq:eucl snmf}
    \Matrix \approx \nonnegFactor \factorCoords
\end{equation}
where $\nonnegFactor \in \Real_{\geq 0}^{\numData \times \numFactors}$ and $\factorCoords \in \Real^{\numFactors \times \dimInd}$. In particular, suppose that the $\sumIndB$th and $\sumIndC$th rows of $\factorCoords$, denoted by $\factorCoords_{\sumIndB,:}$ and $\factorCoords_{\sumIndC,:}$, satisfy $\factorCoords_{\sumIndB,:} = - \factorCoords_{\sumIndC, :}$. Then, if the $\sumIndA$th row of $\nonnegFactor$, corresponding to the $\sumIndA$th element of the dataset, has nonzero $\sumIndB$th and $\sumIndC$th entries, then the $\sumIndA$th element of the dataset could be equally well-represented by a semi-NMF
\begin{equation}
    \Matrix \approx \nonnegFactor' \factorCoords
\end{equation}
where $\nonnegFactor'_{i,j} + \nonnegFactor'_{i,k} = \nonnegFactor_{i,j} + \nonnegFactor_{i,k}$. That is, the approximation can no longer be interpreted as a parts-based decomposition due to the ambiguity about the relationship between the $\sumIndA$th element of the dataset and the $\sumIndB$th and $\sumIndC$th learned factors. Hence, this cancellation effect hinders the interpretability of the representation of the dataset learned by semi-NMF. 

In the case of manifold-valued data $\{\Tensor^\sumIndA\}_{\sumIndA=1}^\numData \subset \manifold$ and their CC-NMDF representation $\log_\mPoint \Tensor^\sumIndA \approx \sum_{\sumIndC=1}^\numFactors \nonnegFactor_{\sumIndA, \sumIndC} \mTVectorFactor_\mPoint^\sumIndC$, an analogous cancellation effect can occur when $(\mTVectorFactor_\mPoint^\sumIndB, \mTVectorFactor_\mPoint^\sumIndC)_\mPoint < 0$ for any pair of $\sumIndB$ and $\sumIndC$. In \Cref{sec:tsnmdf}, we proposed interpreting the factors $\TensorY = \exp_\mPoint \left ( \left (\max_{\sumIndA = 1, \dots, \numData} \nonnegFactor_{\sumIndA,\sumIndC}\right) \mTVectorFactor_\mPoint ^\sumIndC\right )$ as vertices of a nonlinear polyhedron in $\manifold$ that contains the original data points $\Tensor^\sumIndA$. Therefore, we argue that if $(\mTVectorFactor_\mPoint^\sumIndB, \mTVectorFactor_\mPoint^\sumIndC)_\mPoint < 0$ and the learned representation of $\log_\mPoint \Tensor^\sumIndA$ has nonzero coefficients associated with $\mTVectorFactor_\mPoint^\sumIndB$ and $\mTVectorFactor_\mPoint^\sumIndC$ (that is, the $\sumIndA$th row of $\nonnegFactor$ has nonzero $\sumIndB$ and $\sumIndC$ entries), then there exists an ambiguity in the relationship of $\Tensor^\sumIndA$ to the manifold-valued factors $\TensorY^\sumIndB$ and $\TensorY^\sumIndC$. In particular, the vertices associated to the factors may ``cancel out'' with each other, since $(\mTVectorFactor_\mPoint^\sumIndB, \mTVectorFactor_\mPoint^\sumIndC)_\mPoint < 0$. The most extreme example of this occurs when $\mTVectorFactor_\mPoint^\sumIndB = c \mTVectorFactor_\mPoint^\sumIndC$ where $c < 0$. In this case, one of the learned factors is completely redundant, which clearly impacts the effectiveness of the low-rank approximation and the interpretability of the factors. A lower-grade cancellation effect can occur whenever $(\mTVectorFactor_\mPoint^\sumIndB, \mTVectorFactor_\mPoint^\sumIndC)_\mPoint < 0$.

To mitigate this cancellation, we define $\rho:\Real \to \Real$ as 
\begin{equation}
    \rho(x) := \min\{0,x\}.
\end{equation}
Then, we define the \emph{effective} $\sumIndC$th tangent vector factor coordinate of the $\sumIndA$th data point by
\begin{align}
    \nonnegFactorCor_{\sumIndA,\sumIndC} &\coloneqq \nonnegFactor_{\sumIndA,\sumIndC} \mTVectorFactor_\mPoint^\sumIndC + \sum_{\sumIndB\neq \sumIndC}^\numFactors \rho \left( \left( \nonnegFactor_{\sumIndA,\sumIndB} \mTVectorFactor_\mPoint^\sumIndB , \frac{1}{\|\mTVectorFactor_\mPoint^\sumIndC\|_{\mPoint}} \mTVectorFactor_\mPoint^\sumIndC\right)_{\mPoint} \right) \frac{1}{\|\mTVectorFactor_\mPoint^\sumIndC\|_{\mPoint}} \mTVectorFactor_\mPoint^\sumIndC \nonumber \\
    &= \left( \nonnegFactor_{\sumIndA,\sumIndC} + \sum_{\sumIndB\neq \sumIndC}^\numFactors \nonnegFactor_{\sumIndA,\sumIndB} \frac{\rho \left( ( \mTVectorFactor_\mPoint^\sumIndB, \mTVectorFactor_\mPoint^\sumIndC )_{\mPoint} \right)}{\|\mTVectorFactor_\mPoint^\sumIndC\|_{\mPoint}^2} \right) \mTVectorFactor_\mPoint^\sumIndC.
\end{align}
So, using this adjustment, we define corrected manifold-valued factors as
\begin{equation}
    \TensorY^\sumIndC := \exp_{\mPoint}\left( (\max_{\sumIndA=1, \ldots, \numData} \nonnegFactorCor_{\sumIndA,\sumIndC}) \mTVectorFactor_\mPoint^\sumIndC \right).
\end{equation}
We apply this correction throughout our numerical experiments.

\subsubsection{Choice of Base Point of Linearization.} We can also mitigate potential cancellation effects via judicious choice of $\mPoint$, the base point of the chosen tangent space. Since CC-NMDF is a tangent space-based method, the choice of tangent space (via choice of base point $\mPoint$) can have a significant impact on many aspects of the method's performance. For any tangent space-based low rank approximation scheme, a poor choice of base point can result in lower-quality  approximations and identified by examining the exact reconstruction error along with the resulting factors. Such factors may exhibit obvious distortions from expected behavior. In the case of CC-NMDF, the factors may also suffer from cancellation effects, as described in the previous section, if the base point is not chosen carefully. 

Some tangent space-based low-rank approximations of manifold-valued data have natural choices of base points based on their Euclidean analogues; for instance, CC-SVD uses the barycenter of the data as the base point $\mPoint$ because of its relationship to Euclidean principal component analysis \citep{diepeveen2023curvature}. While there is no such analogous choice for the NMF-inspired CC-NMDF, we propose a base point motivated by the previous discussion of cancellation among the learned factors. In particular, one way to ensure that the tangent vector factors learned by CC-NMDF satisfy $(\mTVectorFactor_\mPoint^\sumIndC, \mTVectorFactor_\mPoint^\sumIndB)_\mPoint \geq 0$ is to choose a base point $\mPoint \in \manifold$ such that $(\log_\mPoint \Tensor^\sumIndC, \log_\mPoint \Tensor^\sumIndB)_\mPoint \geq 0$ for all $\sumIndC, \sumIndB \in [\numData]$, such as a point that is sufficiently distant from all of the points in the dataset. However, if such a point is used, then curvature effects can be exaggerated due to the distance between the base point and the data points. As illustrated by \eqref{eq:beta def} and \eqref{eq:relaxed-semi-nmf}, curvature effects are particularly significant for manifolds with negative curvature, since tangent space errors can be amplified on these manifolds. Therefore, in such settings, curvature correction is even more important when using base points satisfying the above heuristic.

\section{Numerics.}
\label{sec:numerics}
To characterize the behavior and applications of the CC-NMDF proposed in \Cref{sec:ccnmdf}, we evaluate our methods on real-world data. Our goal is to compare our scheme to other low-rank approximations for manifold-valued data on the basis of approximation error and interpretability of the resulting factors. We first investigate the effects of different choices of base point and the factor correction proposed in \Cref{sec:algorithm-considerations}. Then, we compare the performance of CC-NMDF to other low-rank approximation schemes for manifold-valued data on the basis of reconstruction quality and interpretability of manifold-valued factors. The code for our algorithms and experiments is available at \url{https://github.com/joycechew/NMDF}.

\subsection{Overview of the Data Set and Experiment Commonalities.} 
\begin{figure}[t!]
    \centering
    \includegraphics[width=0.6\linewidth]{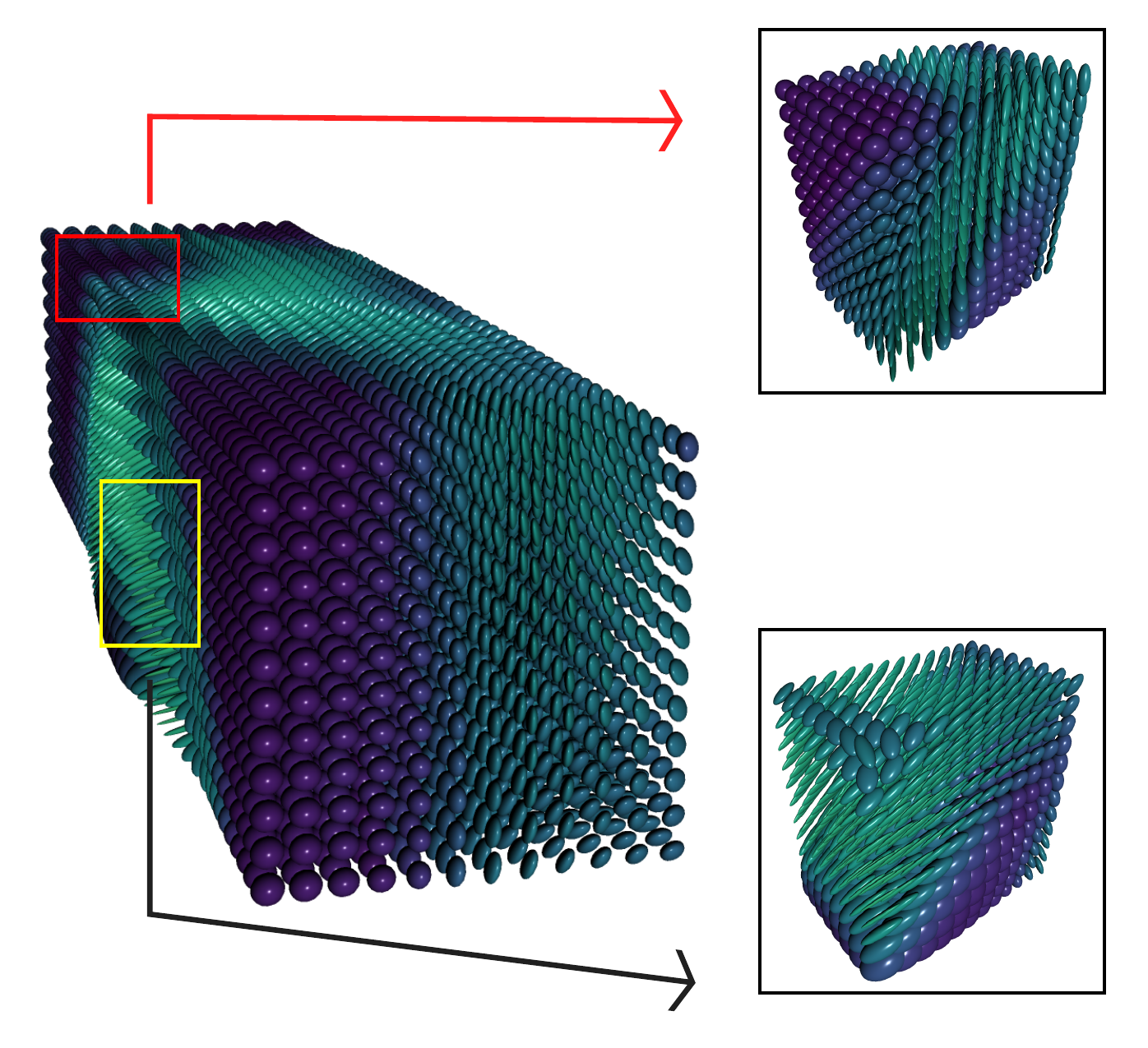}
    \caption{Visualization of the DTI data set (left) and two ``slices'' from the data (right).}
    \label{fig:IIT_dataset}
\end{figure}

We use a data set of diffusion tensor magnetic resonance imaging (DTI) of the adult human brain \citep{zhang2018evaluation,qi2021regionconnect}. In these experiments, we take a subset of the entire brain to make the computation feasible, visualized in \Cref{fig:IIT_dataset}. Each voxel is a $3 \times 3$ symmetric positive definite matrix, so the data lie on $\mathcal{P}(3)$, which is a symmetric Riemannian manifold with non-positive curvature. We focus on a non-compact manifold to investigate our proposed heuristic for choosing a base point. On a compact manifold such as a sphere, it may not be possible to satisfy this heuristic and therefore may become more difficult to identify a good base point. The second reason we focus on data drawn from $\mathcal{P}(3)$ is that we expect from the results of \cite{diepeveen2023curvature} that curvature correction is most notable for non-positively curved spaces, and we want to investigate curvature effects for the proposed factorizations. 

Since the data are arranged as a single mode-3 tensor, we extract and unfold $4 \times 4 \times 4$ regions of voxels and consider each unfolded region an individual data point. Hence, in our processed data set, each data point is on the power manifold $\mathcal{M} = \mathcal{P}(3)^{64}$, and we have $\numData = 117$. We choose to arrange the data in this way to explore the effect of low-rank approximations since this choice of $\mathcal{M}$ has dimension $384$, as opposed to the manifold $\mathcal{P}(3)$, which has dimension $6$. Also, processing the data this way preserves some of the spatial relationships between adjacent voxels in each region, allowing for study of spatial patterns in the data set. 

For each computation of CC-NMDF, we apply the correction described in \Cref{sec:algorithm-considerations} for rendering manifold-valued factors. We perform decompositions for 12 values of $\numFactors$ linearly spaced in $[2,35]$. We run 50 iterations of each iterative method (T-NMDF and CC-NMDF) with $\texttt{maxSubIter} = 5$ for CC-NMDF, and we take $\delta = 0.1$.

Throughout our experiments, we employ the following error metric to evaluate the quality of different low-rank approximations of the data set. First, we define the data power manifold $\manifoldB = \manifold^{N}$ and use the distance inherited from the power manifold structure
\begin{equation}
    \distance_{\manifoldB}(\Tensor, \TensorY) = \sqrt{\sum_{\sumIndA=1}^\numData \distance_{\manifold}(\Tensor^\sumIndA, \TensorY^\sumIndA)^2}
\end{equation}
as our error metric. That is, given data $\{\Tensor^\sumIndA\}_{\sumIndA=1}^\numData \subset \manifold$ and tangent space approximations $\Xi_\mPoint^\sumIndA \approx \log_\mPoint \Tensor^\sumIndA$, we take the approximation error to be
\begin{equation}
    \distance_{\manifoldB} (\Tensor, \Xi_\mPoint) = \sqrt{\sum_{\sumIndA=1}^\numData \distance_{\manifold}(\Tensor^\sumIndA, \exp_{\mPoint}(\Xi_{\mPoint}^\sumIndA))^2}.
\end{equation}

\subsection{Choice of Base Point for Linearization.}
\label{sec:numerics base point}
We first explore the effect of the base point for linearization on the resulting CC-NMDF. A natural first choice for the base point would be the barycenter of the data, since this choice is successfully used for SVD-like approximations of manifold-valued data \citep{diepeveen2023curvature}. However, as discussed in \Cref{sec:algorithm-considerations}, this choice may result in redundant factors or otherwise negatively impact the quality of the factors. Another possible choice of base point is a point $\mPoint$ on the manifold sufficiently distant from each data point such that $(\log_{\mPoint} \Tensor^\sumIndA, \log_{\mPoint} \Tensor^\sumIndB)_{\mPoint} > 0$ for all $\sumIndA \neq \sumIndB$. Often, one suitable choice of this kind of base point is a point close to zero on the manifold. In this case, we use the point on $\mathcal{P}(3)^{64}$ with every entry equal to $10^{-5} \times I_3$ as such a point, where $I_3$ is the $3 \times 3$ identity matrix, and we denote this point by $\mPointB$. For comparison, we let $\mPointC$ denote the barycenter of the data, and we compute the CC-NMDF of our data set using $\mPointB$ and $\mPointC$ as the base points of linearization. We denote these as CC-NMDF($\mPointB$) and CC-NMDF($\mPointC$), respectively, in our results. 

\begin{figure}[t!]
    \centering
    \includegraphics[width=0.48\linewidth]{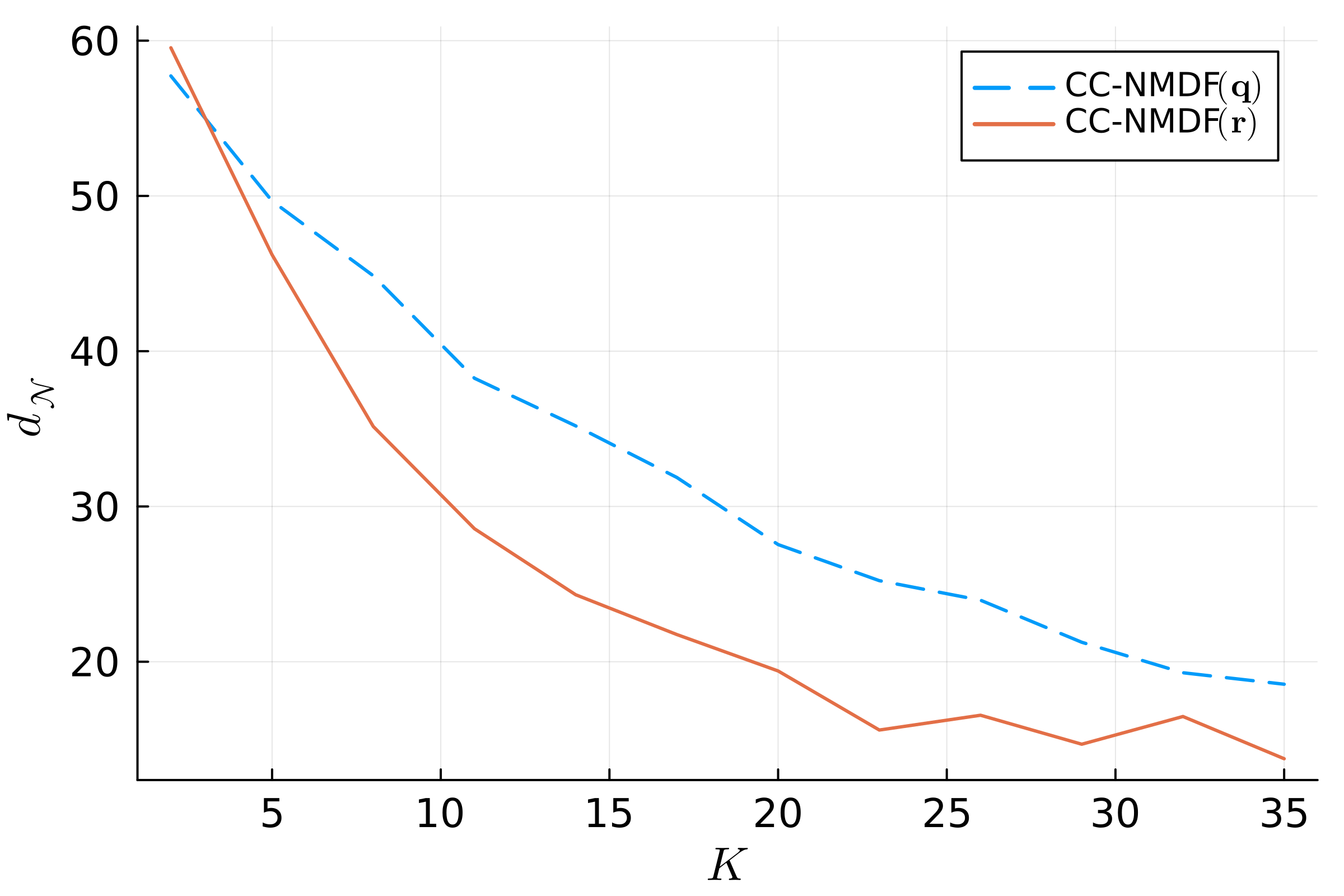}
    \caption{The error incurred by CC-NMDF using approximate zero ($\mPointB$) and the data barycenter ($\mPointC$) as the base point of linearization, plotted against approximation rank ($\numFactors$).}
    \label{fig:cc_q0_qmean_loss}
\end{figure}
\begin{figure}[t!]
    \centering
    \begin{subfigure}{0.48\linewidth}
        \includegraphics[width=\linewidth]{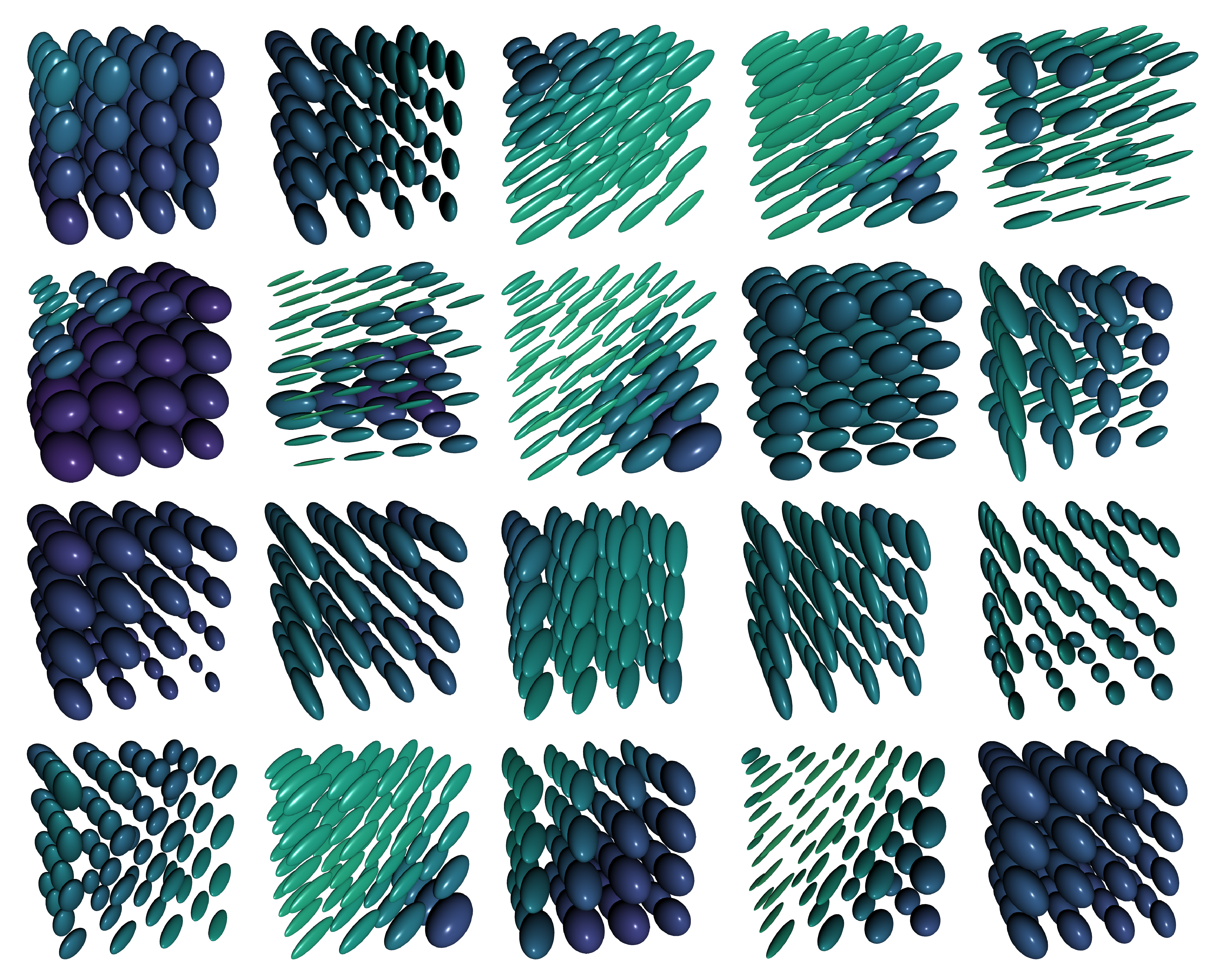}
        \caption{}
    \end{subfigure}
    \hfill
    \begin{subfigure}{0.48\linewidth}
        \includegraphics[width=\linewidth]{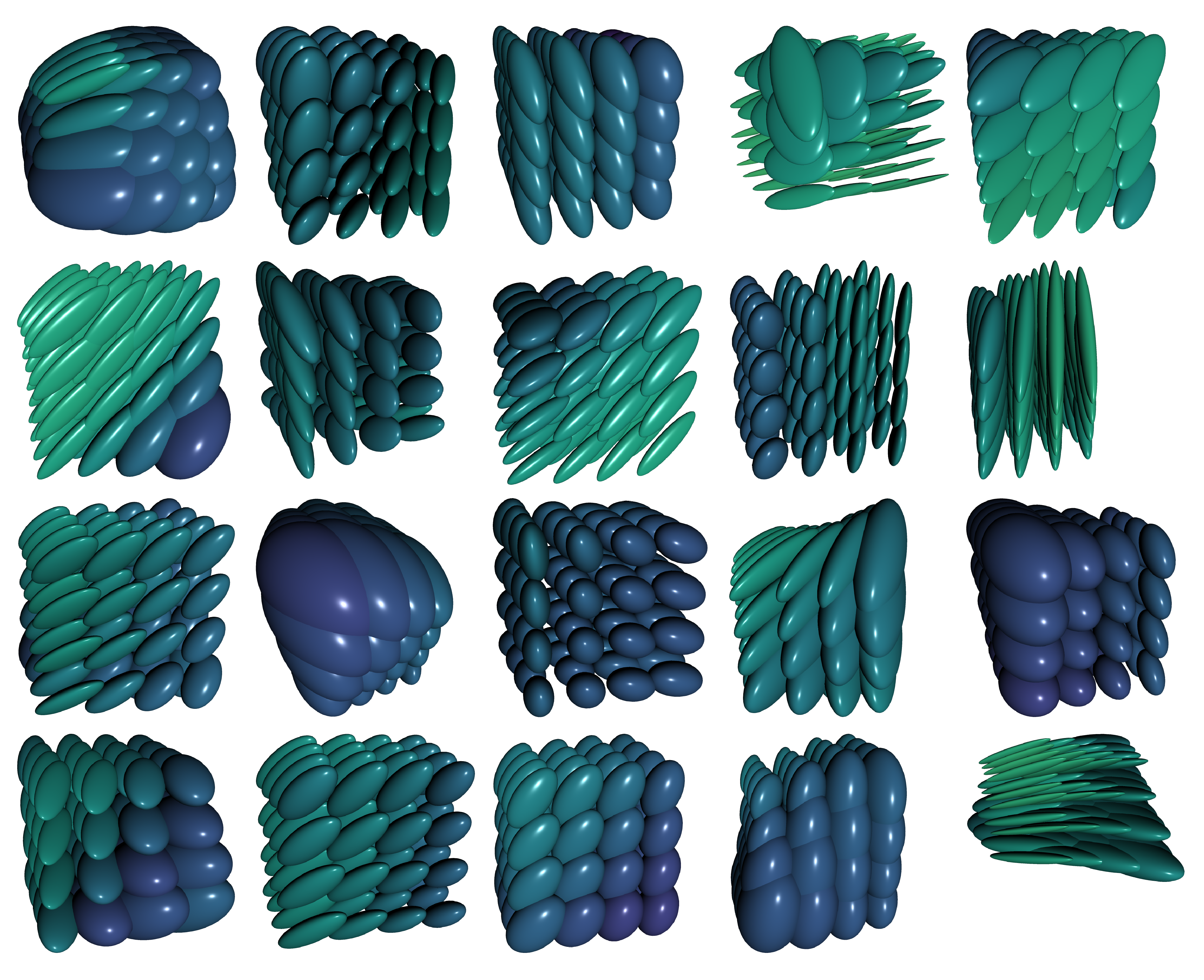}
        \caption{}
        \label{fig:cc_qmean_factors}
    \end{subfigure}
    \caption{The factors obtained from rank-20 CC-NMDF using (a) $\mPointB$, approximately zero, and (b) $\mPointC$, the data barycenter, as the base point of linearization. (b) shows a clear deterioration in the interpretability of the factors when compared to (a) and the original data set.}
    \label{fig:cc_q0_qmean_factors}
\end{figure}

We see that using $\mPointC$ as the base point of linearization results in a better approximation error (\Cref{fig:cc_q0_qmean_loss}), but using $\mPointB$ results in much more interpretable manifold-valued factors (\Cref{fig:cc_q0_qmean_factors}). The choice of base point $\mPoint$ impacts the effect of distortion due to the geometry of $\manifold$ on the linearized data $\{\log_\mPoint \Tensor^\sumIndA\}_{\sumIndA=1}^\numData$ that is used in the computation of low-rank approximations. A base point that is further from the manifold data points will naturally induce more distortion in the linearization itself, so we expect that CC-NMDF($\mPointC$) will have a better approximation error that CC-NMDF($\mPointB$). However, the manifold-valued factors arising from CC-NMDF($\mPointC$) as the base point of linearization are subject to distortions that do not appear to be representative of actual features in the original data set, and the individual factors appear to generally be more uniform and lack detail that is present in the manifold-valued factors from CC-NMDF($\mPointB$). Hence, on the basis of interpretability of the manifold-valued factors, we use $\mPoint$ as the base point of linearization for manifold NMF methods in the following section.

\subsection{Comparison with T-NMDF and CC-SVD.}
\label{sec:numerics vs naive}
\begin{figure}[t!]
    \centering
    \begin{subfigure}{0.48\linewidth}
        \includegraphics[width=\linewidth]{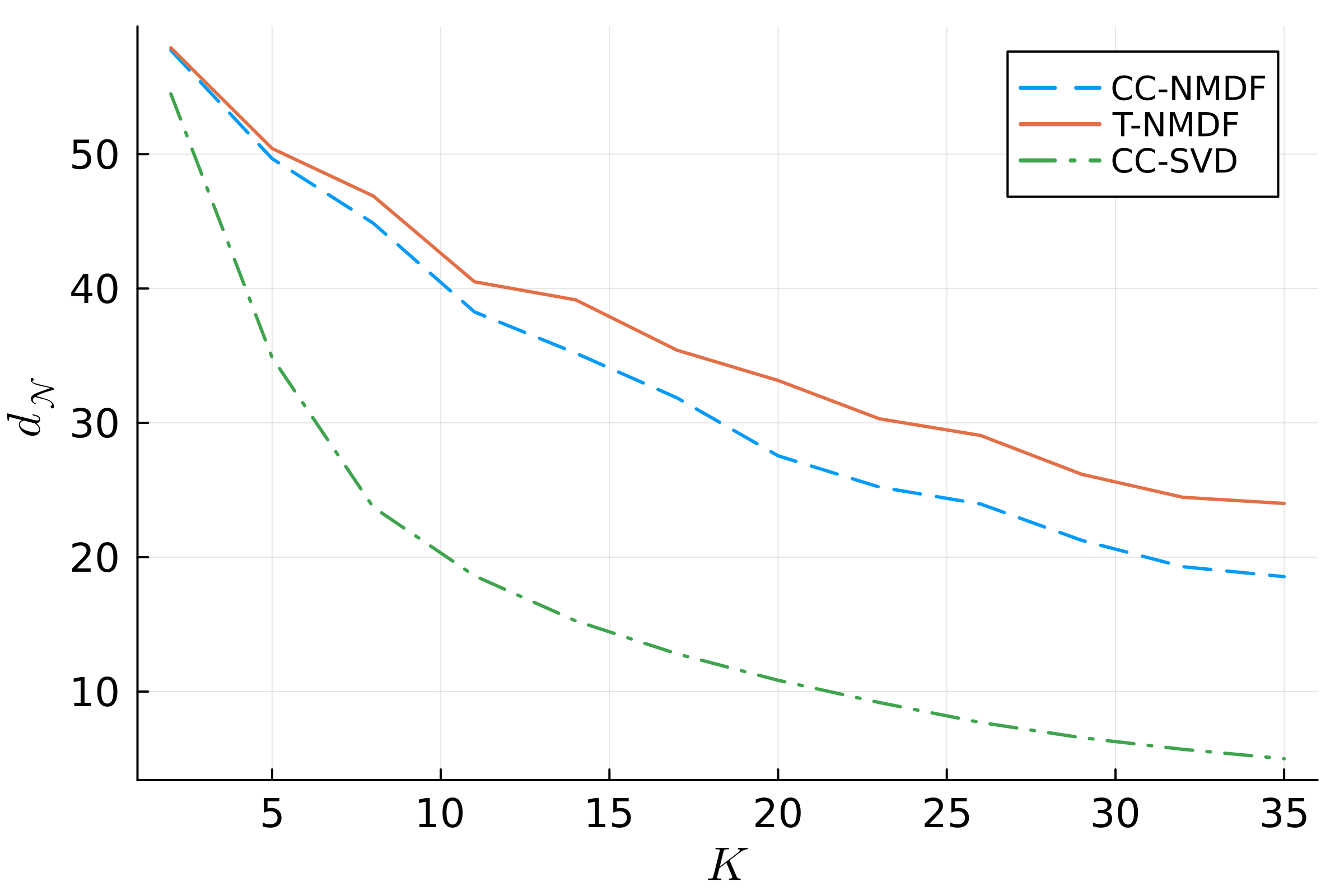}
    \caption{}
    \label{fig:naive_vs_cc_loss}
    \end{subfigure}
    \hfill
    \begin{subfigure}{0.48\linewidth}
        \includegraphics[width=\linewidth]{figures/cc_q0_factors.png}
        \caption{}
        \label{fig:cc_q0_factors}
    \end{subfigure}
    \begin{subfigure}{0.48\linewidth}
        \includegraphics[width=\linewidth]{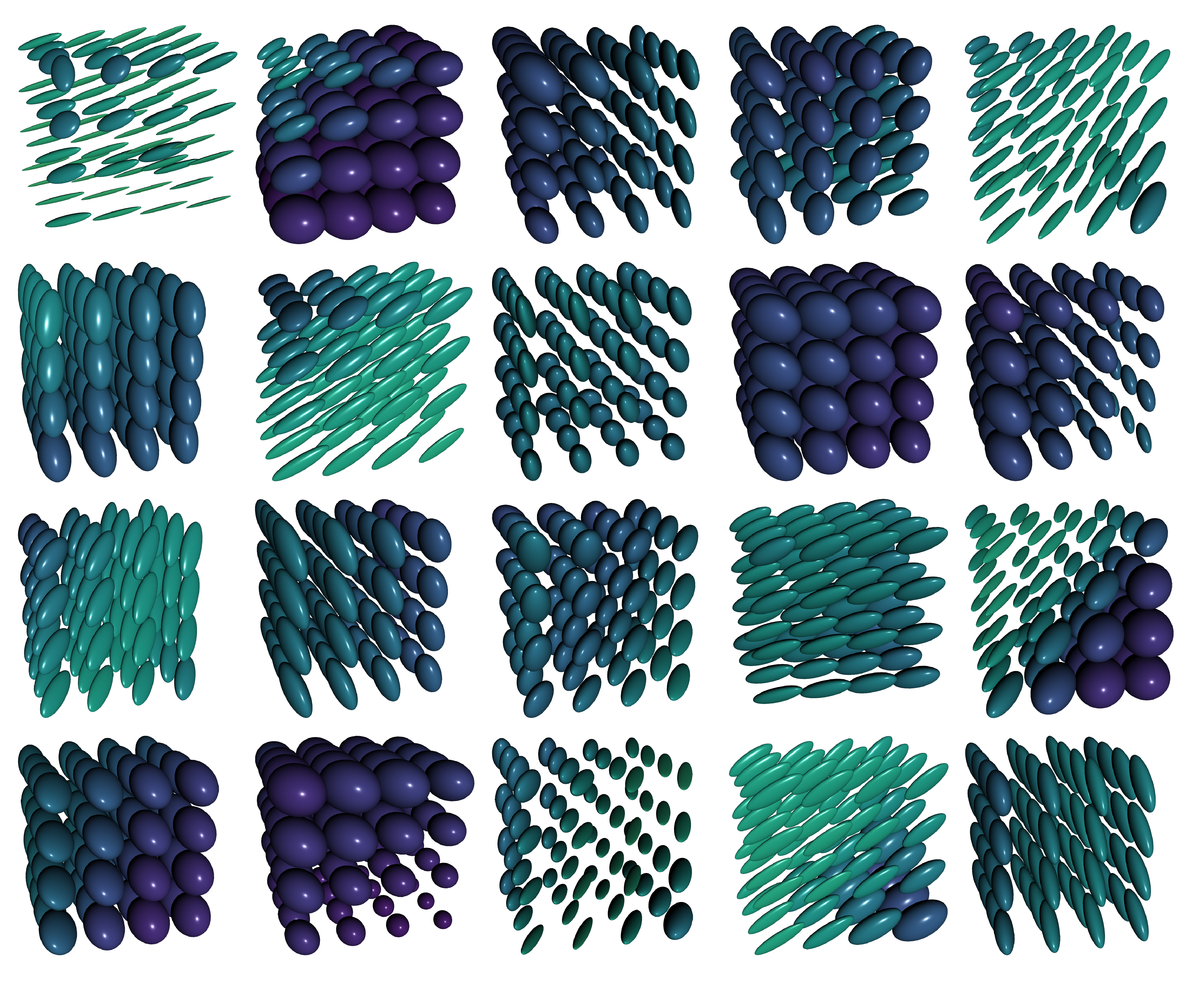}
        \caption{}
        \label{fig:naive_q0_factors}
    \end{subfigure}
    \hfill
    \begin{subfigure}{0.48\linewidth}
        \includegraphics[width=\linewidth]{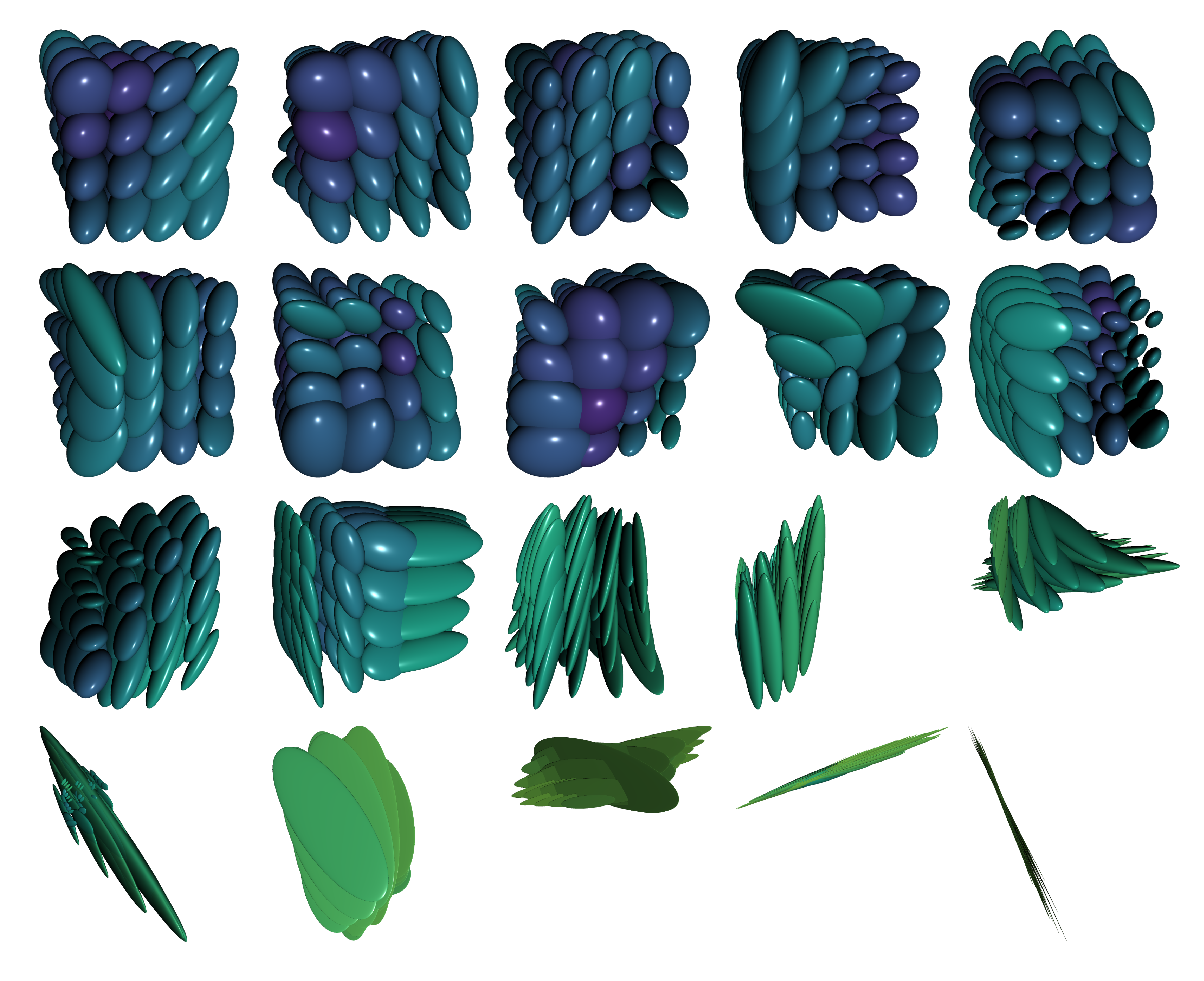}
        \caption{}
        \label{fig:CC-SVD_qmean_factors}
    \end{subfigure}
    \caption{(a) Approximation error incurred by CC-NMDF, T-NMDF, and CC-SVD for approximation ranks 2-35. (b) Manifold-valued factors obtained by CC-NMDF. (c) Manifold-valued factors obtained by T-NMDF. (d) Manifold-valued factors obtained by CC-SVD.}
\end{figure}
Next, we compare CC-NMDF to T-NMDF, the uncorrected tangent space-based method described in \Cref{alg:t-SNMDF}, and CC-SVD, a curvature-corrected singular value decomposition for manifold-valued data \citep{diepeveen2023curvature}. T-NMDF is simpler to implement and requires less memory and computational time than CC-NMDF, as it does not require any manifold mappings or expansion into a basis of $\tangent_\mPoint \manifold$. However, as a purely tangent space-based method, T-NMDF does not account for the curvature of the underlying manifold.

In \Cref{fig:naive_vs_cc_loss}, we see that the approximation error on the manifold is moderately better for CC-NMDF compared to T-NMDF. Since we expect non-negligible curvature effects on $\mathcal{P}(3)^{64}$, it is unsurprising that the curvature-corrected method results in a better approximation error. Furthermore, several factors that result from T-NMDF appear to be similar to each other (see \Cref{fig:naive_q0_factors}), indicating that T-NMDF does not capture as many features of the data set as CC-NMDF does. Additionally, some of the factors appear to suffer from distortions that do not correspond to features in the original data set, such as the factors in the upper left and upper right of \Cref{fig:naive_q0_factors}. Therefore, the curvature correction incorporated in CC-NMDF improves both the manifold approximation error and interpretability of the resulting factors. 

It is unsurprising that CC-SVD has the best approximation error on the manifold, as it is designed to emulate the behavior of the Euclidean SVD, and it is the solution of an unconstrained low-rank approximation problem. However, in \Cref{fig:CC-SVD_qmean_factors}, we see that the factors obtained via CC-SVD (\Cref{fig:CC-SVD_qmean_factors}) sharply contrast with the factors obtained via CC-NMDF (\Cref{fig:cc_q0_factors}). In general, the CC-SVD factors do not correspond well to any of the individual regions in the original data set. Furthermore, the factors with the smallest singular values (shown in the bottom two rows of \Cref{fig:CC-SVD_qmean_factors}) are extremely distorted, indicating that they capture either noise or other pathologies in the data. 

\section{Conclusions.}
\label{sec:conclusions}
We have introduced an analogue of nonnegative matrix factorization for manifold-valued data, which we refer to as CC-NMDF. Our CC-NMDF incorporates the geometry of the underlying manifold domain to learn manifold-valued factors that faithfully represent characteristic features of the data in the original domain. At the same time, our method operates on a tangent space of the underlying manifold, so it only requires the computation of relatively few manifold mappings. We provide an iterative algorithm to compute CC-NMDF, which alternates between solving a linear system and a multiplicative update. In our experiments, we find that CC-NMDF produces manifold-valued factors that are interpretable as characteristic features of the input dataset. Our work opens the door to adapting variants of NMF to manifold-valued data, enabling more powerful and interpretable processing of these data that remain faithful to the geometry of the data. Future work on manifold NMF-type methods should also investigate how to choose an appropriate base point for general manifolds, since data on compact manifolds may not admit a base point that satisfies our proposed heuristic.

\acks{WD was partially supported by the UCLA Dunn Family Endowed Chair funds.  DN was partially supported by NSF DMS 2408912 and JC was partially supported by NSF DMS 2136090 and NSF DGE 2034835.}

\bibliography{Bibliography}

\end{document}